\documentstyle[12pt,amsfonts,amscd]{amsart}
\newtheorem{Theorem}{Theorem}[section]
\newtheorem{Definition}[Theorem]{Definition}
\newtheorem{Proposition}[Theorem]{Proposition}

\newtheorem{Lemma}[Theorem]{Lemma}
\newtheorem{Corollary}[Theorem]{Corollary}
\theoremstyle{remark}
\newtheorem{Example}[Theorem]{Example}

\def\il{\int\limits_}
\def\CC{{\Bbb C}}
\def\RR{{\Bbb R}}
\def\ZZ{{\Bbb Z}}

\def\l{\lambda}
\def\LL{\Lambda}

\def\ovr{\overline}

\def\Dl{\Delta}
\def\dl{\delta}

\def\bd{\partial}

\def\sbs{\subset}

\def\Aut{\operatorname{Aut}}

\def\be{\begin{enumerate}}
\def\ee{\end{enumerate}}
\def\bT{\begin{Theorem}}
\def\eT{\end{Theorem}}
\def\bP{\begin{Proposition}}
\def\eP{\end{Proposition}}
\def\bD{\begin{Definition}}
\def\eD{\end{Definition}}
\def\bE{\begin{Example}}
\def\eE{\end{Example}}
\def\bL{\begin{Lemma}}
\def\eL{\end{Lemma}}
\def\bC{\begin{Corollary}}
\def\eC{\end{Corollary}}

\def\H{{\mathcal H}}

\begin{document}
\title{Perturbation of domains and automorphism groups}
\author{Buma L. Fridman and Daowei Ma}
\begin{abstract} The paper is devoted to the description of changes of the structure of the holomorphic automorphism group of a bounded domain in ${\Bbb C}^n$ under small perturbation of this domain in the Hausdorff metric. We consider a number of examples when an arbitrary small perturbation can lead to a domain with a larger group, present theorems concerning upper semicontinuity property of some invariants of automorphism groups. We also prove that the dimension of an abelian subgroup of the automorphism group of a bounded domain in ${\Bbb C}^n$ does not exceed $n$.
\end{abstract}

\address{buma.fridman@@wichita.edu, Department of Mathematics,
Wichita State University, Wichita, KS 67260-0033, USA}
\address{dma@@math.twsu.edu, Department of Mathematics,
Wichita State University, Wichita, KS 67260-0033, USA}
\maketitle \setcounter{section}{-1}
\section{Introduction}Let $D$ be a bounded domain in $\CC^n$, and $\Aut(D)$ the automorphism group of $D$ (the group of biholomorphic self-maps of $D$). It is well known that  $\Aut(D)$ is a (real) Lie group. In this article we examine the following question. If a domain $\tilde D$ is a small perturbation of
$D$, what is the relation between $\Aut(D)$ and $\Aut(\tilde D)$? The automorphism group of a domain is its group of symmetries, and a natural idea (coming from Euclidean geometry) is that ``small perturbations may destroy symmetry, but never create symmetry''. So, one expects that if a perturbation is small enough then automorphism groups will satisfy an upper semicontinuity principle: $\Aut(\tilde D)$ is isomorphic to a subgroup of  $\Aut(D)$. In particular one would not expect an ``increase'' of the $\Aut(D)$ by arbitrarily small perturbations of $D$. 

In the early eighties R.~Greene and S.G.~Krantz ([GK1], [GK2], [GK3]) examined this question. They proved an upper semicontinuity result in the $C^2$ topology, and gave the first counterexample to the upper semicontinuity principle. In [Ma], [FP], [FMP] the semicontinuity question has been examined further for various other topologies. We also note here that in Riemannian geometry, results of this kind have been obtained by R. Palais [Pa], D. Ebin [Eb] and K. Grove-H. Karcher [GKa]. 

In this paper we summarize some of the results concerning upper semicontinuity of the automorphism groups in the topology induced by the Hausdorff metric on bounded domains in $\CC^n$, and present new theorems and examples complementing those results. The paper consists of two parts. In the first part we collect counterexamples to the upper semicontinuity principle. In the second part we present theorems proving the upper semicontinuity property for some invariants of the automorphism groups, and state applications of these results.

\setcounter{section}{0}
\section{Violations of the upper semicontinuity principle}
\subsection{A counterexample in the topology induced by the Hausdorff metric} Greene-Krantz in [GK3] proved the following upper semi-continuity theorem for the $C^2$ topology on domains in $\CC^n$ (see Theorem 2 in that paper).
\bT\label{T:mt} If $D$ is a $C^2$ strongly pseudoconvex domain in $\CC^n$ that is not biholomorphic to the ball, then there exists a neighborhood of $D$ in the $C^2$ topology, such that for any $C^2$ domain  $\tilde D$ in that neighborhood $\Aut(\tilde D)$ is isomorphic to a subgroup of $\Aut(D)$.\eT

Below we present an example that shows that in the Hausdorff metric the above theorem will not hold. The example is in $\CC$. (${\Bbb Z}_k$ denotes the cyclic group containing k elements).
\bE
There is a bounded $C^\infty$ domain $D$ in $\CC$ that is not simply connected, and such that for any neighborhood $U$ of $D$ in the Hausdorff metric there is a $C^\infty$ domain $\tilde D$ in $U$ such that $\Aut(\tilde D)$ is not isomorphic to any subgroup of $\Aut(D)$.
\eE

\noindent {\it Construction.} $B(z_0,r)$ denotes an open disk with center $z_0$ and radius $r.$ Consider an $(N+1)$-connected domain $D=\Delta \setminus \bigcup\limits_{s=1}^N\overline{\Delta}_s$ where $\Delta =B(0,1)$ is the open unit disk and all $\overline{\Delta}_s$ are smaller non-intersecting closed disks, whose boundaries lie entirely in $\Delta$. For a given $1>\varepsilon >0$ fix a positive $ \varepsilon _1\leq \varepsilon $ and such that the set $S=\{z\in \Delta \mid Re(z)>-1+\varepsilon _1\}$ contains all  $\overline{\Delta}_s$. Suppose a natural number $j>1$ is also given. We now choose a positive $\delta $ such that $L(S)\subset B(1,1/2^j)$ where $L$ is a M\"{o}bius transformation $L(z)=\frac{z+a}{1+za}$ , and $a=1-\delta $. We observe that $L(\Delta )=\Delta $.  Consider now $M=$ $\bigcap\limits_{k=0}^{j-1}(L(D)\cdot\exp(\frac{2\pi k}ji))$ (each term is a rotation of $L(D)$ by angle $\frac{2\pi k}j$). One can verify that by construction $\ZZ_j$ acts on $M$. We define $\widetilde{D}_j=L^{-1}(M)$. Then $\widetilde{D}_j\subset D$, and since $\Aut(D_j)\cong \Aut(M)$, $\ZZ_j$ is isomorphic to a subgroup of $\Aut(D_j)$. Also the difference $D\setminus \widetilde{D}_j\subset \Delta \setminus S$ and therefore the Hausdorff distance between $D$ and $\widetilde{D}_j$ is less than $\varepsilon $. If $N\geq 2$ the group $\Aut (D)$ is finite. Since $j$ could be chosen to be arbitrarily large, the statement has been proved.

\noindent {\it Remarks}. The above construction works for any finitely connected domain: in any neighborhood of this domain and any integer $j$ there is a domain whose automorphism group contains ${\Bbb Z}_j$. (One could compare this with Theorem 1.4 below). A similar construction can be done in $\CC^n$ for any $n\geq 1$.

As mentioned in the introduction, in [GK2] the authors provide a counterexample of the failure of upper semicontinuity in the $C^{1-\varepsilon}$ topology (for any $\varepsilon >0$):
\bE
There are pseudoconvex domains $\{D_j\}_{j=1}^\infty $ and $D_0$, each of which is $C^\infty $ and strongly pseudoconvex except at one point, such that $\Aut(D_j)\neq \{id\}$ for all $j\ge 1$, $\Aut(D_0)=\{id\}$, and $D_j\longrightarrow
D_0$ in the $C^{1-\varepsilon}$ topology, any $\varepsilon>0 $.
\eE

\subsection{Groups that act on a dense set of domains} Let $\H^n$ be the metric space of all bounded domains in ${\Bbb C}^n$ with the metric equal to the Hausdorff distance between boundaries of domains. In [FP] the following phenomenon is described: groups $\ZZ_k$ act on a dense set of domains in $\H^n$. We present this in the following statement (follows from [FP], Theorem 2.1):
\bT
For any positive integers $n$ and $k$, the set of bounded domains in  $\CC^n$ such that their automorphism group contains a subgroup isomorphic to ${\Bbb Z}_k$ is everywhere dense in $\H^n$.
\eT

\noindent {\it Remark}. By examining the proof of Theorem 2.1 in [FP], one can prove the above statement for the infinite discrete group $\ZZ$: for any positive integer $n$ the set of bounded domains in  $\CC^n$ such that their automorphism group contains a subgroup isomorphic to $\ZZ$ is everywhere dense in $\H^n$. We will leave it to the interested reader to fill in the necessary details in the proof of Theorem 2.1 in [FP] and obtain the proof of this statement. 

So, for any domain (even a rigid one, i.e. with $\Aut(D)=\{id\}$) one can make a perturbation of less than a given size and obtain a domain with a large, even an infinite, cyclic group. The natural question arises: which Lie groups will a similar statement hold for? Below in this section we provide a positive answer for any finite group, and a negative answer for connected groups of positive dimension.

For any finite group the following statement holds. Let $G$ be a group of order $m<\infty $ .

\bT For any $n\geq m$ the set of bounded domains in   $\CC^n$ whose automorphism group contains a subgroup isomorphic to $G$ is everywhere dense in $\H^n$.
\eT
\begin{pf}
1. Group G, as a group of order $m$, is isomorphic to a subgroup of the permutation group $P_m$ of $m$ elements. It is therefore sufficient to prove the statement for the space $\CC^n$ and the group $P_m$ , which acts naturally on the unit ball of $\CC^n$ (and on the entire $\CC^n$) by simply permuting the first $m$ coordinates.

2. Any $\H^n$ neighborhood of a given bounded domain in $\CC^n$ contains a bounded domain which is a finite union of open balls. This follows from the definition of a domain as an open connected set.

3. Let $D$ be a domain which is a finite union of open balls, and $U$ a neighborhood of $D$. Consider the smallest ball $B_1$ in $\CC^n$ that contains $D$. By using first a linear transformation $L$ we can map $B_1$ onto the unit ball $B=B(0,1)$ and then by using an element $M\in \Aut(B)$ we can assure that the points $p,q$: $p=(1,0,...,0),q=-p$, both belong to the boundary of the image $(M\circ L)(D)$. Now by adding to this image small portions of the unit ball near $p,q$ we proceed to construct a domain $D_1\subset B$, such that (1) $D_1\cap (B(p,\delta)\cup B(q,\delta))=B\cap (B(p,\delta)\cup B(q,\delta))$ for some positive $\delta $, and (2) for the biholomorphic transformation $F=L^{-1}\circ M^{-1}$ the following holds: $F(D_1)\subset U$.

4. All that's left to prove is that the theorem holds for a domain $D_1$ with the above property: $D_1$ is a subset of $B$ that contains full open $\delta $-portions of $B$ near the points $p,q$. Let $U_1$ be an $\varepsilon >0$ neighborhood of $D_1$ in the metric space $\H^n$, we may assume $\varepsilon <\min (\delta, 1/8)$. Now we choose a positive $\alpha<1$ such that $T(S)\subset B(p,\varepsilon)$, where $S=\{z\in B\mid |z-q|\geq \varepsilon /2\}$, $T$ is the following automorphism of the unit ball $
\{z_1{\rightarrow}\frac{z_1+a}{1+z_1a},z_k\rightarrow \sqrt{1-a^2}\frac{z_k}{1+z_1a},k=2,...,n\}$, $a=1-\alpha $. Define now $D_2=\bigcap\limits_{\sigma \in P_m}\sigma (T(D_1))$, $\widetilde{D}_1=T^{-1}(D_2)$. Note that by construction $P_m$ is a subgroup of $\Aut(D_2)\cong\Aut(\widetilde{D}_1)$. One can now also check that by construction $\widetilde{D}_1\subset U_1$.
\end{pf}

To improve the above theorem we now discuss the following. If $G$ is a subgroup of $\Aut(D)$, can one make a perturbation of $D$ within a given neighborhood, so that the resulting domain has automorphism group {\it isomorphic} to $G$? The following lemma gives a partial answer to this question for a finite subgroup. 

\bL\label{L:id} Let $D$ be a bounded domain in $\CC^n$ without isolated boundary points. Let $G$ be a finite subgroup of $\Aut(D)$. Then for each neighborhood $U$ of $D$ in $\H^n$ there is a domain $\widetilde{D}\in U$, such that $\Aut(\widetilde{D})\cong G$.
\eL

\begin{pf} For a set $X\subset D$, the orbit of $X$ with respect to $G$ is denoted $O_G(X)=\bigcup\limits_{g\in G}g(X)$. For an automorphism $g\in \Aut(D)$, we denote by $S_g=\{x\in D| g(x)=x\}$ the set of fixed points for this automorphism. For $g\neq id$ the set $S_g$ is closed and nowhere dense in $D$. We choose now a point $x\in D\backslash \bigcup\limits_{g\in G,g\neq id}S_g$. Then $g_1(x)\neq g_2(x)$, for any two different $g_1,g_2\in G$. If the point $x$ is chosen close enough to the boundary of $D$, then the domain $\widehat{D}=D\backslash O_G(x)$ is such that $\widehat{D}\in U$; we choose $x$ to satisfy this condition also. Denote by $\Theta (x,r)$ the ball in the Bergman metric on $D$ with center $x$ and radius $r.$ Consider now a small $\delta >0$, and $n+1$ points $x_j\in \Theta (x,\delta )$ different from $x$ such that the following four conditions hold:

(1) all the pairwise distances (in the Bergman metric) between different images $g(\Theta (x,\delta)),g\in G$ are greater than $\delta$, 

(2) $[D\backslash O_G(\Theta (x,\delta ))]\subset U$,

(3) the points $\{x_j\}_{j=1}^{n+1}$ are in ``general position'', so that if any automorphism of $D$ fixes these points, then this automorphism is the identity (see [FKK]), and 

(4) all pairwise Bergman distances for the $n+2$ points $x,\{x_j\}_{j=1}^{n+1}$ are different. 

Consider now $\widetilde{D}=[D\backslash O_G(x\cup (\bigcup\limits_{j=1}^{n+1}x_j))].$ Because of the properties listed above, $\widetilde{D}\subset U,$ and $G$ acts on $\widetilde{D}$. Suppose $h\in \Aut(\widetilde{D})$. Because of the definition of $\widetilde{D}$, $h$ can be automatically extended to an automorphism of $D$, so $h\in \Aut(D).$ It is also clear that (see properties (1), (4) above) $h(x)\in O_G(x).$ So there is a $g\in G$, such that $g(x)=h(x)$. We have $(h^{-1}\circ g)(x)=x$ and, because of the positions of points $x_j$, $(h^{-1}\circ g)(x_j)=x_j$.  It follows from Property (3) that $h^{-1}\circ g=id$. So, $h=g$. We have proved that $\Aut(\widetilde{D})\cong G$.
\end{pf}

By examining the proof of Theorem 1.5, one can see that the domains constructed there had no isolated boundary points, and $P_m$ was a subgroup of their automorphism group. Therefore, by using the above lemma we get the following refinement of the theorem 1.5:
\bT For any $n\geq m$ the set of bounded domains in   $\CC^n$ whose automorphism group is isomorphic to $G$, is everywhere dense in $\H^n$.
\eT

As far as which Lie groups can be realized as automorphism groups of bounded domains: in [SZ], [BD] it is proved that any compact Lie group can be realized as the group of automorphisms of a smooth strictly pseudoconvex domain, and in [TS] it is shown that any linear Lie group can be realized as the group of automorphisms of a bounded domain.

Theorems 1.4-1.7 show that arbitrarily small perturbation of a domain in $\H^n$ may create a domain with a larger automorphism group. But, for all known examples, this group is discrete, so it is of dimension zero. The natural question arises: can small perturbation in $\H^n$ create domains with larger \textit{dimensions} of automorphism groups? Though we have a general theorem 2.1 (see also Corollary 2.2) stated later it is reasonable to give an example (using an argument completely different from the proof of the theorem 2.1) that shows that domains whose automorphism groups are connected and of positive dimension are not everywhere dense in $\H^n$. 

\bT
If $D$ is a strongly pseudoconvex domain in $\CC^n$ with a discrete automorphism group, then there exists a neighborhood of $D$ in the $\H^n$ topology, such that for any domain  $\tilde D$ in that neighborhood $\Aut (\tilde D)$ is not a connected group of positive dimension.
\eT

In the proof we will use a non-negative invariant function $h_Q(z)$ defined on a bounded domain $Q$ and introduced in [Fr1]. The function has the following properties.

\begin{itemize}
\item $h_Q(z)$ is invariant under biholomorphic transformations.
\item $h_Q(z)$ is continuous.
\item If $h_Q(z_0)=0$ for some $z_0\in Q$, then $Q$ is biholomorphic to the unit ball of $\CC^n$.
\item For a strictly pseudoconvex domain $D$, $h_D(z)\to 0$ as $z$ approaches the boundary.
\item If $D_j\to D$ in the Hausdorff metric, then $h_{D_j}\to h_D$ uniformly on compacta.
\end{itemize}

\begin{pf}
Let $\l=\max_{z\in D} h_D(z)$,  $z_0\in D$ be such that $h_D(z_0)=\l$, $W=\{z\in D: h_D(z)\ge \l/2\}$,  $K=\{z\in D: h_D(z)= \l/2\}$, let $U_+=\{z\in D: h_D(z)>\l/2\}$, and  $U_- =\{z\in D: h_D(z)< \l/2\}$. 

If the statement does not hold, then there is a sequence $\{D_j\}$ converging to $D$ in the topology induced by the Hausdorff distance such that $\Aut(D_j)$ is connected and has positive dimension. Since $h_{D_j}\to h_D$ and $h_{D_j}$ is constant on $X_j=\Aut(D_j)(z_0)$, the orbit $X_j$ does not intersect $K$ for large $j$. It follows that $X_j\subset U_+\cup U_-$. Since $X_j$ is connected, $X_j\cap U_+$ or $X_j\cap U_-$ is empty. But $z_0$ is in the former, so the latter must be empty. Therefore, $X_j\subset W$ for large $j$. By Corollary~4.1 in [FP], $\Aut(D_j)$ is isomorphic to subgroup $\Aut(D)$, which is a contradiction. 
\end{pf}

\subsection{A different counterexample} A different kind of a counterexample to the upper semicontinuity principle is exhibited below: $\Aut(D)$ is isomorphic to $\RR$, but the automorphism groups of some nearby domains are isomorphic to $S^1$. This example appears in [FMP], we describe it here in full for completeness.

\bE
There is a sequence $\{D_j\}$ of bounded pseudoconvex domains in $\CC^2$ converging in $\H^2$ to a domain $D$ such that $\Aut(D_j)\cong  S^1$ for each $j$, and $\Aut(D)\cong\RR$.
\eE

{\it Construction.} Let $\Dl$ denote the unit disc in $\CC$. Let $Q_j=\{z\in \Dl: |z- 2^{-1}+2^{-j}|>1/2\}$, $Q=\{z\in \Dl: |z- 2^{-1}|>1/2\}$, $D_j=\{(z,w): z\in Q_j, w\in \Dl, w\ne z\}$,  $D=\{(z,w): z\in Q, w\in \Dl, w\ne z\}$.

1. One can see that $D_j\to D$.

2. The domains $D_j$ and $D$ are bounded and pseudoconvex.

3. We now prove that $\Aut(D)\cong\RR$. Let $F\in \Aut(D)$. On each fiber $(z,\cdot)$, $F$ is bounded and has an isolated singularity, so $F$ extends to be an automorphism of $Q\times \Dl$. Thus, $F$ has the form $F(z,w)=(f(z), g(w))$, or $F(z,w)=(g(w), f(z))$. For both cases, one has, by the definition of $D$, that $$(*)\;\;\;\;\;\;\;\;\;\; f(z)=g(z),\;\;\;  z\in Q$$. The second case is impossible, since $f(Q)= \Dl$, $g(\Dl)=Q$, and $f(Q)=g(Q)$ with lead to a contradiction that $\Dl$ coincides with a subset of $Q$. Therefore, $F$ has the form $F(z,w)=(f(z), g(w))$, where $f\in\Aut(Q)$, $g\in\Aut(\Dl)$. By (*), $f=g|_Q$. Let $\phi(w)=-i(w+1)/(w-1)$. Then $\phi$ is a biholomorphic map from $\Dl$ to the upper half-plane $\Pi=\{\zeta\in \CC: Im(\zeta)>0\}$, and $\phi(Q)=\LL\equiv\{\zeta\in \CC: 0<Im(\zeta)<1\}$. Now $\phi\circ g\circ \phi^{-1}$ is an automorphism of $\Pi$, and its restriction to $\LL$ is an automorphism of $\LL$. Thus $\phi\circ g\circ \phi^{-1}(\zeta)=\zeta +t$ for some $t\in \RR$. It follows that $\Aut(D)=\{F_t: t\in \RR\}\cong \RR$, where $F_t(z, w)=(g_t(z), g_t(w))$, and
$$g_t(w)=\phi^{-1}(\phi(w)+t)={2w+i(w-1)t\over 2+i(w-1)t}.$$

4. In a way very similar to the above argument, one can prove that $\Aut(D_j)\cong  S^1$ for each $j$. 
\endpf

\section{Upper semicontinuity results in $\H^n$ and applications}
\subsection{The upper semicontinuity result for dimensions of automorphism groups}
In addition to the upper semicontinuity result of Greene-Krantz in the smooth topology in $\CC^n$ stated in the beginning of Section 1, there exist analogous positive statements in [FP], [Ma] confirming the upper simicontinuity principle for domains in topologies different from the smooth topology; we will not restate them here. Instead we will focus on results in the $\H^n$ topology. As pointed out in Section 1 the upper semicontinuity principle for the automorphism groups does not hold in $\H^n$, but the principle will be true for the {\it dimensions} of the groups. The following is the main semicontinuity result in [FMP].

\bT\label{T:mt} The function $\dim\Aut(D)$ is upper semicontinuous on $\H^n$.
\eT

The reader is referred to [FMP] for the proof of this theorem. An immediate consequence is the following

\bC\label{C:gc} For
each $k>0$ the set of all domains in $\H^n$ whose groups of
automorphisms have dimensions greater than or equal to $k$ is closed
and, therefore, nowhere dense.\eC

Thus a domain cannot be approximated by domains whose automorphism groups have strictly larger dimensions. In the next section we present other consequences of Theorem~2.1

\subsection{Approximation of domains and automorphism groups}

We say that a domain $D$ can be approximated by a domain $M$, if it can be approximated by biholomorphic images of $M$ in the metric space $\H^n$, e.g. if any neighborhood of $D$ (in $\H^n$) contains a biholomorphic image of $M$. We call $U$ a universal domain if any bounded domain in $\CC^n$ can be approximated by $U$; the existence of such domains was established in [Fr2]. The following statements are direct consequences of the Theorem 2.1.

\bC\label{C:gc1}
$\dim\Aut(U)=0$ for any universal domain $U$.\eC

\bC\label{C:gc2}
If a domain {D} can be approximated by a domain {M} then $\dim(\Aut(M))\le\dim(\Aut(D))$.\eC

\bC\label{C:gc3} If each of the two domains $D_1$ and $D_2$ can be approximated by the other then $\dim(\Aut(D_1))=\dim(\Aut(D_2))$.\eC

To complement this statement we provide the following
\bE 
There exist two pseudoconvex domains $D_1$, $D_2$ in $\CC^n$, such that each can be approximated by the other (in $\H^n$) but their automorphism groups are not isomorphic.
\eE

\noindent {\it Construction.}We will extensively use a universal domain $D$ construction provided in [Fr2].  By examining that construction we note that $D$ may be chosen to satisfy all three properties listed below: 

(1) $D$ does not have isolated boundary points.

(2) $D$ can be chosen as $D=B\backslash S$, where $B=B(0,1)$ and $S$ is a closed set, $S\subset B(p,1/4)$, $p=(1,0,...,0)$ (see Lemma 2.1 in [FP] to verify this property). 

(3) If $M\subset \{z | Re(z_1)<1/2\}$, and $D\backslash M$ is a domain, then $D\backslash M$ is also a universal domain.

We first prove our statement for $n=1$, so we start with such a universal domain $D\subset \CC$. 

1. Fix three points $\zeta _1,\zeta _2,\zeta _3\in B(0,1/2)$, such that all three pairwise Kobayashi distances in $D$ between them are different. Define $\check{D}_1=D\backslash (\zeta _1\cup \zeta _2\cup \zeta _3)$. Since any automorphism of $\check{D}_1$ extends to $\zeta _j,$ the choice of these points leads to $\Aut(\check{D}_1)=\{id\}$. (see also [PL])

2. Define $\check{D}_2=D\cap (-D).$ By construction $\ZZ_2$ acts on $\check{D}_2,$ so $\Aut(\check{D}_2)\neq \{id\}$. 

3. Since both $\check{D}_1,\check{D}_2$ are universal domains in $\CC$, they approximate each other, and the statement is proved for $n=1$.

4. For any $n>1$, consider $D_k=\check{D}_k\times B^{n-1}$, for $k=1,2$ where $B^{n-1}
$ is the unit ball in $\CC^{n-1}$. One can now check that the statement will hold for $D_1,D_2$. This finishes our construction.

The automorphism group of any Reinhardt domain has at least dimension n, and therefore the class of biholomorphic images of Reinhardt domains is nowhere dense in $\H^n$; moreover

\bC\label{C:gc4} If a domain $D\subset\subset \CC^n$ and $\dim(\Aut(D))<n$ then there is a neighborhood of $D$ that contains no biholomorphic images of Reinhardt domains.\eC

Analogously, for circular domains the following statement holds:

\bC\label{C:gc5}
If a domain $D\subset\subset \CC^n$ and $\dim(\Aut(D))=0$ then there is a neighborhood of $D$ that contains no biholomorphic images of circular domains.\eC

\subsection{On abelian subgroups of automorphism groups} We start with the following

{\bf Definition.} The {\it rank} of a bounded domain $D$ in ${\Bbb C}^n$ is the maximum of the dimensions of the closed abelian subgroups of $\Aut (D)$.

Let $r(D)$ denote the rank of $D$.

\bT The function $r(D)$ is upper semicontinuous on $\H^n$.
\eT

\begin{pf}
Let $D_j$ be a sequence of domains converging in $\H^n$ to a domain $D$. Let  us choose a ball $\hat B=B(p,r+a)$, $r,a>0$, such that the closure of $\hat B$ belongs to all $D_j$ for sufficiently large $j$. Let $B=B(p, r)$. We may assume that the ranks of all $D_j$ are the same and equal to $k$. Choose vector fields $X^m_j$, $1\le m\le k$, in the Lie algebra of $\Aut(D_j)$ such that 
$$[X^m_j,X^l_j]=0;\;\;\;\il{\hat B}(X^m_j,\ovr X^l_j)\,dV=\dl_{ml},$$ 
where $\delta_{ml}$ is Kronecker's delta.
 
By Lemma~2.4 in [FMP], $\|X_j^m\|_{\hat B}\le C\|X_j^m\|_{B}$ for some positive constant $C$ independent of $j, m$. On the other hand, for some positive constant $c=c(n)$, $\|X_j^m\|_B\le c a^{-n}$, since for each $z\in B$,
$$|X_j^m(z)|\le v_{n}^{-1}a^{-2n}\int_{B(z,a)} |X_j^m|\,dV$$
$$\le v_n^{-1/2}a^{-n}\Big(\int_{\hat B} |X_j^m|^2\,dV\Big)^{1/2}=v_n^{-1/2}a^{-n},$$
where $v_n$ is the volume of the unit ball in ${\Bbb C}^n$. It follows that for some positive constants $C_1, C_2$, 
$$1\ge C_1\|X_j^m\|_{B}\ge C_2\|X_j^m\|_{\hat B}.$$

Let $g^m_j$ be the one-parameter groups generated by $X^m_j$. By Theorem~2.5 in [FMP] one can choose a subsequence  $\{j_k\}$ such that $g^m_{j_k}$ converge, uniformly on compacta in $D\times \RR$, to a one-parameter group $g^m(z,t)$, and $X^m_{j_k}$ converge to a vector field $X^m$ uniformly on compacta in $D$. Since
$$[X^m, X^l]=0, \;\;\;\il {\hat B}(X^m,\ovr X^l)\,dV=\dl_{ml},$$
the closed group generated by $g^m(z,t)$, $m=1,\dots, k$, is an abelian subgroup of $\Aut(D)$ of dimension at least $k$. Therefore, $r(D)\ge k$.
\end{pf}

The maximum dimension of $\Aut(D)$ for a domain $D\subset\subset \CC^n$ is $n^2+2n$. However the dimension of the abelian subgroup of any such automorphism group is much less, e.g. the following exact estimate holds.
\bT
If $D$ is a bounded domain in ${\Bbb C}^n$, then $r(D)\le n$.
\eT

\begin{pf}
Seeking for a contradiction suppose that $r(D)=k>n$. Let $H$ be a connected $k$-dimensional abelian subgroup of $\Aut(D)$. Then $H$ is isomorphic to $T^m\times {\Bbb R}^{k-m}$, $0\le m\le k$. Choose vector fields $X_1(z), \dots, X_k(z)$ in the Lie algebra of $H$ so that $X_1(z), \dots, X_m(z)$ generate a subgroup $Q$ of $H$ isomorphic to $T^m$, and $X_{m+1}(z), \dots, X_k(z)$ generate a subgroup of $H$ isomorphic to ${\Bbb R}^{k-m}$. By the Slice Theorem in [Br], there exists a $z_0\in D$ such that $Qz_0$ is diffeomorphic to $T^m$. Thus $X_1(z_0), \dots, X_m(z_0)$ are linearly independent over $\Bbb R$. Since $k>n$, there are complex constants $a_j=b_j+ic_j$, $j=1,\dots, k$, not all 0, such that $\sum_{j=1}^k a_jX_j(z_0)=0$. Let $U=\sum_{j=1}^k b_jX_j$, $V=\sum_{j=1}^k c_jX_j$. Then $V(z_0)=iU(z_0)$.

We now show that $U(z_0)\ne 0$. Suppose that $U(z_0)=0$. Then $V(z_0)=0$. 
Therefore, $z_0$ is a fixed point of the one-parameter subgroup generated by $U$. Thus that one-parameter subgroup belongs to a compact subgroup of $H$, hence $b_{m+1}=\cdots=b_k=0$. Since $X_1(z_0), \dots, X_m(z_0)$ are linearly independent over $\Bbb R$, we see that $b_1=\cdots=b_m=0$. Thus, $b_1=\dots=b_k=0$. Similarly, from $V(z_0)=0$ we derive $c_1=\cdots=c_k=0$. It follows that $a_1=\cdots=a_k=0$, contradicting the choice of $a_j$. Therefore, $U(z_0)\ne 0$.

Let $g_u(z,t)$, $g_v(z,t)$ be the one-parameter groups generated by $U$, $V$ respectively. Define a map $f: {\Bbb C}\to D$ by $f(t+is)=g_v(g_u(z_0,t),s)$. Since $[U, V]=0, \;[U,iU]=0$, we see that 
$$V(f(t+is))=\frac{\partial g_v}{\partial z}(g_u(z_0, t),s)\frac{\partial g_u}{\partial z}(z_0,t)V(z_0),$$
$$iU(f(t+is))=\frac{\partial g_v}{\partial z}(g_u(z_0, t),s)\frac{\partial g_u}{\partial z}(z_0,t)iU(z_0).$$
The above two equations, together with $V(z_0)=iU(z_0)$, yield that $V=iU$ on the image $f(\Bbb C)$. It follows that $f$ is a nonconstant holomorphic map from $\Bbb C$ to $D$, contradicting Liouville's Theorem. Therefore, $r(D)\le n$.
\end{pf}

\subsection{Is the characteristic number of $\Aut(D)$ upper semicontinuous on $\H^n$?} 
By Iwasawa's theorem (see [MZ, p.~188]) the identity component of the group $\Aut(D)$ is homeomorphic to $K\times{\Bbb R}^{k}$, where $K$ is a maximal compact subgroup and $k$ is the {\it characteristic number} of $\Aut(D)$. It is interesting to find out what happens with $K$ and ${\Bbb R}^k$ under small perturbations of domains. It seems to us that the characteristic number of $\Aut(D)$ is also semicontinuous in $\H^n$. Denote this number by $char(\Aut(D))$. Below we make two remarks related to this problem.

At this time we can prove the following statement, that follows from Theorem 3.2 in [FMP]. We include the proof here for completeness.
\bP Let $D\sbs \CC^n$ be a bounded domain. If each neighborhood of $D$ contains a $\tilde D$ such that the $char(\Aut(\tilde D))\ge 1$, then the $char(\Aut(D))\ge 1$.
\eP

\noindent {\it Remark}. This statement shows that $S^1$ and $\RR$ cannot be interchanged in example 1.9. 

\begin{pf}
By the hypothesis there is a sequence $\{D_j\}$ of domains converging to $D$ such that for each $j$, the identity component $G_j$ of $\Aut(D_j)$ is noncompact. Seeking for a contradiction, suppose that the identity component $G$ of $\Aut(D)$ is compact. Fix a $z_0\in D$. The orbit $G(z_0)$ is compact. We may assume that $G(z_0)\sbs D_j$ for each $j$. For each connected component $H$ of $\Aut(D)$, either the set $H(z_0)$ coincides with $G(z_0)$ or $G(z_0)\cap H(z_0)=\emptyset$. Indeed, if $h\in H$ and $h(z_0)\in G(z_0)$, then $H(z_0)=Gh(z_0)=G(z_0)$, since $H=Gh$. Now we claim that there exists a positive number $a$ such that $a< d(H(z_0), G(z_0))$ for each component $H$ of $\Aut(D)$ with $H(z_0)\ne G(z_0)$. Otherwise, there is a sequence $\{H_k\}$ of distinct components of $\Aut(D)$ with $H_k(z_0)\ne G(z_0)$ such that $d(H_k(z_0), G(z_0))\to 0$. Passing to a subsequence if necessary, we may assume that there are $h_k\in H_k$ such that  $h_k(z_0)$ tends to a point in $G(z_0)$. It follows that some subsequence of $\{h_k\}$ converges in the compact-open topology to a $g\in \Aut(D)$; but this is impossible because $h_k$ belong to different components of the Lie group $\Aut(D)$. Therefore, such an $a$ exists. Decreasing $a$ if necessary, we see that the open set
$$V=\{z\in D: d(z, G(z_0))<a\}$$
is relatively compact in $D$ and each $D_j$, and satisfies $\overline V\cap \Aut(D)(z_0)=G(z_0)$. This implies that $\bd V\cap \Aut(D)(z_0)=\emptyset$. Since $G_j$ is noncompact, $G_j(z_0)$ is noncompact, hence $G_j(z_0)\cap \bd V\ne \emptyset$. It follows that for each $j$ there is a $g_j\in G_j$ with $g_j(z_0)\in \bd V$. Some subsequence of the sequence $\{g_j\}$ converges uniformly on compacta to a $g\in \Aut(D)$. It is clear that $g(z_0)\in \bd V$, contradicting $\bd V\cap \Aut(D)(z_0)=\emptyset$. Therefore, $G(D)$ is noncompact.
\end{pf}

The example below (in $\CC^3$) shows that the proof of the possible semicontinuity theorem for the characteristic number cannot be expected to be done in a straightforward manner. As proved above, if a bounded domain $D$ is a limit of bounded domains $D_j$ whose automorphism groups have closed one-parameter subgroups, that are non-compact (e.g. $\cong\RR$) then the $\Aut(D)$ also has such a subgroup. The example below shows however that it is not necessary for any such subgroup of $\Aut(D)$ to be a limit of any subsequence of the corresponding subgroups of $\Aut(D_j)$.

\bE\label{E:1} There is a sequence $\{D_j\}$ of bounded domains in $\CC^3$ converging (in  $\H^3$) to a domain $D$ such that
$\Aut(D_j)\cong\RR$ for each $j$, and $\Aut(D_j)\to G$, $G$ is a one-dimensional subgroup of $\Aut(D)$ and $G\cong\ S^1$.
\eE

\noindent {\it Construction}. Let $A=\{z_1:1/2<|z_1|<1\},$ $B=\{z_2:||Im(z_2)|<1\},Q_j=Q_j(z_1,z_2)=\{z_3:z_3\in A,z_3\neq p_j,q_j\},$ where $p_j=p_j(z_1,z_2)=(z_1/(|z_1|+1/(j+1))/(\exp (ijz_2)/|\exp (ijz_2)|)$, $q_j=1-1/(j+1)^2.$ $D_j=\{(z_1,z_2,z_3):z_1\in A,z_2\in B,z_3\in Q_j(z_1,z_2)\}$, $D=A\times B\times A$.

1. One can see that $D_j\to D$.

2. We now describe $\Aut(D_j).$ First consider $g_j(z,t)=(\exp(it)z_1,z_2+t/j,z_3).$ One can verify that $g_j(z,t)\in Aut(D_j)$ for each $t\in\RR.$ Also, the group $G_j=\{g_j(z,t),t\in\RR \}\cong\RR.$ We now prove that $\Aut(D_j)=G_j.$ We have already verified $\Aut(D_j)\supseteq G_j.$ Pick now $g\in \Aut(D_j).$ By construction one can see that $g(z_1,z_2,z_3)=(f_1(z_1),f_2(z_2),f_3(z_3)),$ where $f_1,f_2$ are automorphisms of $A,B$ respectively, and $f_3$ maps $Q_j(z_1,z_2)$ conformally onto $Q_j(f_1(z_1),f_2(z_2)).$ Since $f_3$ can be extended holomorphically to $p_j(z_1,z_2),q_j(z_1,z_2)$, $f_3$ after this extension becomes an automorphism of $A.$ Considering all the possible choices for $f_3$, we see that this automorphism has to be the identity, and the only choice for $f_3($ $p_j(z_1,z_2))$ is $p_j(f_1(z_1),f_2(z_2))$. So, $p_j(z_1,z_2)=p_j(f_1(z_1),f_2(z_2))$ for all $(z_1,z_2)\in A\times B$. $f_1\in \Aut(A),$ and therefore can be one of the following: (a) $f_1=\exp (i\theta )/2z_1$ ; (b) $f_1=\exp (i\theta )$ for some $\theta\in{\Bbb R}$.

Consider now the equation $p_j(z_1,z_2)=p_j(f_1(z_1),f_2(z_2))$ which has to
hold for all $(z_1,z_2)\in A\times B$.

By using explicit expressions (taking the absolute value of both sides of the equation) one can prove that case (a) cannot take place.

The case (b) leads to $|\exp (ij(z_2-f_2(z_2))|=\exp (i[j(z_2-f_2(z_2))+\theta ])$, and since the last expression is an analytic function it must be a positive constant. So, there is such an integer $k,$ and a real constant $c,$ that $j(z_2-f_2(z_2))+\theta =2\pi k+ic.$ Therefore, $f_2(z_2)=z_2+(\theta -2\pi k)/j-ic/j$. Since $f_2$ is an automorphism of the horizontal strip $B,$ $c=0.$ Consider now $\alpha =\theta -2\pi k,$ then $f_1(z_1)=\exp (i\alpha )z_1,f_2(z_2)=z_2+\alpha /j,$ and, as we noticed before, $f_3(z_3)=z_3$.

We have proved that $g\in G_j.$ So, $\Aut(D_j)\subseteq G_j,$ and therefore $\Aut(D_j)=G_j$.

4. Note now that $\Aut(D_j)=G_j\to G$, where  $G=\{g(z,t)=(\exp (it)z_1,z_2,z_3),t\in{\Bbb R}\}$, and therefore $G\cong\ S^1$.

5. The domains $D_j$ and $D$ are not bounded. However if $f$ is a Riemann map from $B$ to the unit disk, then $F=(z_1, f(z_2), z_3)$ maps all $D_j$ and $D$ biholomorphically onto bounded domains, and those domains will have the property that was the purpose of this example.


\begin{thebibliography}{999}
\bibitem[BD]{BD} E. Bedford, J. Dadok, {\sl Bounded domains with
prescribed group of automorphisms,} Comment. Math. Helv. {\bf 62} (1987), 561--572
\bibitem[Br]{Br} G.~Bredon, {\sl Introduction to compact transformation groups}, Academic Press, New York, 1972
\bibitem [Eb]{Eb} D. Ebin, {\sl The manifold of Riemannian metrics}, Global
analysis, Proceedings of Symposium in Pure Mathematics, XV, AMS (1970), 17-40 
\bibitem[FP]{FP} B. L. Fridman, E. A. Poletsky, {\sl Upper
semicontinuity of automorphism groups}, Math. Ann., 299(1994),
615--628
\bibitem[FMP]{FMP} B. L. Fridman, D. Ma, E. A. Poletsky, {\sl Upper
semicontinuity of the dimensions of automorphism groups in $\CC^n$}, to appear in Amer. J. Math 125 (2002)
\bibitem[Fr1]{Fr1} B. L. Fridman, {\sl Biholomorphic invariants of a hyperbolic manifold and some applications}, Trans. Amer. Math. Soc., 276 (1983), no.2, 685-698
\bibitem[Fr2]{Fr2} B. L. Fridman, {\sl A universal exhausting domain}, Proc. Amer. Math. Soc., 98 (1986), 267--270
\bibitem [FKK]{FKK} B. L. Fridman, K. T. Kim, S. G. Krantz, D. Ma, {\sl On fixed points and determining sets for holomorphic automorphisms}, to appear in Michigan Math. J. 50 (2002)
\bibitem[GK1]{GK1} R. Greene, S.G. Krantz, {\sl The Automorphism groups of
strongly pseudoconvex domains}, Math. Ann., 261(1982), 425--446 
\bibitem[GK2]{GK2} R. Greene, S.G. Krantz, {\sl Stability of the Caratheodory and Kobayashi metrics and applications to biholomorphic mappings,} Proceedings of Symposia in Pure Math. Providence: AMS vol.41 (1984), 77-94
\bibitem[GK3]{GK3} R. Greene, S.G. Krantz, {\sl Normal Families and the Semicontinuity of Isometry and Automorphism Groups,} Math. Z., vol.190 (1985), 455-467
\bibitem[GKa]{GKa} K. Grove, H. Karcher, {\it How to conjugate $C^1$-close
group actions}, Math. Z., 132, (1973), 11--20
\bibitem[Ko]{Ko} Sh. Kobayashi, {\sl Hyperbolic Complex Spaces},
Springer-Verlag, 1998
\bibitem[Ma]{Ma} D. Ma, {\sl Upper semicontinuity of isotropy and
automorphism groups}, Math.\ Ann., 292(1992), 533-545
\bibitem[MZ]{MZ} D. Montgomery, L. Zippin, {\sl Topological transformation
groups}, New York, Interscience, 1955
\bibitem[Pa]{Pa} R. Palais, {\sl Equivalence of nearby differentiable actions
of a group}, Bull. Amer. Math. Soc., 67, (1961), 362--364
\bibitem[PL]{PL} E. Peschl and M Lehtinen, {\sl A conformal self-map which
fixes 3 points is the identity}, Ann.\ Acad.\ Sci.\ Fenn., Ser. A I Math.,  4
(1979), no.~1, 85--86
\bibitem[SZ]{SZ} R. Saerens, W. R. Zame, {\sl The isometry groups of
manifolds and the automorphism groups of domains,} Trans. Amer. Math.
Soc. {\bf 301} (1987), 413--429
\bibitem[TS]{TS} A. E. Tumanov, G. B. Shabat, {\sl Realization of
linear Lie groups by biholomorphic automorphisms of bounded domains},
Funct. Anal. Appl., 24 (1990), 255--257
\end{thebibliography}
\end{document}